\documentclass[letterpaper,10pt]{amsart}
\usepackage[bb=fourier]{mathalfa}
\usepackage{amsmath}
\usepackage{amssymb}
\usepackage{mathtools}      % improves amsmath
\usepackage{stmaryrd} % the commands \llbracket, \rrbracket replace \ldbrack, \rdbrack from package mathabx
\usepackage{esint} % extended set of integrals for Computer Modern
\usepackage{enumitem}
\usepackage{amsthm}
\usepackage{thmtools}
\usepackage{tikz}
\usepackage{tikz-cd}
\usepackage[colorlinks=true,backref=page,bookmarksopen=true]{hyperref} % Links in the pdf. Backref option: each bibliographical entry denotes where it was cited
\usepackage[capitalise]{cleveref} % Must be loaded after hyperref. It's advisable to compile the file with LUALATEX instead of PDFLATEX to avoid some annoying warning messages.
\usepackage[font=small]{caption}

\usepackage{iftex}
\ifPDFTeX
   \usepackage[utf8]{inputenc}
   \usepackage[T1]{fontenc}
   \usepackage{lmodern}
\else
   \ifXeTeX
     \usepackage{fontspec}
   \else
     \usepackage{luatextra}
   \fi
   \defaultfontfeatures{Ligatures=TeX}
\fi

\usepackage{xcolor}
\hypersetup{
    colorlinks,
    linkcolor={red!50!black},
    citecolor={blue!50!black},
    urlcolor={blue!80!black}
}

\renewcommand*{\backref}[1]{}
\renewcommand*{\backrefalt}[4]{\quad \tiny 
  \ifcase #1 (\textbf{NOT CITED.})%
  \or    (Cited on page~#2.)%
  \else   (Cited on pages~#2.)%
  \fi}
\begin{document}

\title[Boundary and interior solutions of the cohomological equation]{Relating boundary and interior solutions of the cohomological equation for cocycles by isometries of negatively curved spaces. The Li\v vsic case.}
\date{\today}

\author[A.~Moraga]{Alexis Moraga}
\address{Facultad de Matem\'aticas, Pontificia Universidad Cat\'olica de Chile\\
Av. Vicuna Mackenna 4860, Macul, Chile}
\email{\href{mailto:ajmoraga@mat.uc.cl}{ajmoraga@mat.uc.cl}}

\author[M.~Ponce]{Mario Ponce}
\address{Facultad de Matem\'aticas, Pontificia Universidad Cat\'olica de Chile\\ Av. Vicuna Mackenna 4860, Macul, Chile}
\email{\href{mailto:mponcea@mat.uc.cl}{mponcea@mat.uc.cl}} 


\begin{thanks}
{The authors were partially supported by CONICYT PIA ACT172001 and  by Proyecto FONDECYT 1180922.}
\end{thanks}

\maketitle

\begin{abstract}
We consider the reducibility problem of cocycles by isometries of Gromov hyperbolic metric spaces in the Li\v vsic setting. We show that provided that the boundary cocycle (that acts on a compact space) is reducible in a suitable  H\"older class, then the original cocycle by isometries (that acts on an unbounded space) is also reducible.       
\end{abstract}

\section{Introduction}

The study of the dynamical properties of a cocycle is simplified substantially when we can reduce it to a cocycle that takes values in a simpler group. More precisely, given a homeomorphism (dynamical system) $T : \Omega \to \Omega$ over a compact metric space $\Omega$ and a (topological)
group $\mathcal{G}$, we consider a continuous $\mathcal{G}$-valued cocycle  $A :  \Omega \to
\mathcal{G}$ over the base dynamics $T$.    We are interested in the dynamical behavior of the cocycle, that is, we are interested in the properties of the dynamical product $A^n(\omega):=A(T^{n-1}\omega)\cdot A(T^{n-2}\omega)\cdots A(\omega)$. Such kind of cocycles (and its dynamical properties) appear in many situations in dynamical systems an other branches of mathematics. For instance, when $T$ is a diffeomorphism of a smooth $n$-manifold, the cocycle given by the derivative $A(\omega):=DT(\omega)\in GL(n, \mathbb{R})$ provides a lot of information about the dynamics of $T$. In fact, this is one the main approaches to the theory of smooth dynamical systems. The source of examples for the study of these cocycles is vast, including the Schr\"odinger equation. In that case, the underlying group is $SL(2, \mathbb{R})$ and the (base) dynamical system $T$ corresponds to a minimal linear translation on a torus (see  \cite{AVKR06}, \cite{FK}). The Kontsevitch-Zorich cocycle is a keystone for studying the Teichmuller flow over translation surfaces (see \cite{YOCC}, ). \\

To reduce a cocycle consists in to find a continuous function $B:\Omega\to \mathcal{G}$ such that the conjugated cocycle $B(T\omega)^{-1}\cdot A(\omega)\cdot B(\omega)$ takes values in a small subgroup of $\mathcal{G}$. In the case of the trivial subgroup $\{e_{\mathcal{G}}\}$, where $e_{\mathcal{G}}$ is the neutral element, we look for the existence of a function $B:\Omega\to \mathcal{G}$ such that 
\begin{equation}\label{eq3}
B(T\omega)\cdot B(\omega)^{-1}=A(\omega).
\end{equation}
We call this equation the {\it cohomological equation}. In that case  we say that $A$ is a {\it coboundary} or that $A$ is  {\it reducible}.  Note that (\ref{eq3}) yields $A^n(\omega)=B(T^{n-1}\omega)\cdot B(\omega)^{-1}$. \\

The challenges about reducibility of cocycles are usually divided according to the characteristics of the base dynamics $T$. Namely, there exist good results regarding if $T$ is elliptic, hyperbolic or a partially hyperbolic dynamical system. For instance, in the minimal case (let us think to rotations), a lot of work had been assembled under the Theory of Gottschalk and Hedlund (see \cite{GOHE}, \cite{KARO},  \cite{MMY}), and the KAM theory (see \cite{AVKR06}). In the hyperbolic case, the Li\v vsic theory takes account of many results (see  \cite{delallave2010}, \cite{dolgo05}, \cite{forfla}, \cite{KAL10}, \cite{KS}, \cite{KARO}, \cite{LIV72}, \cite{NavasPonce}, \cite{nittor1}, \cite{poll99}).  Recently the solution to the cohomological equation in the abelian case has been addressed by some authors as a central problem in the partially hyperbolic dynamical system case (see \cite{KAKO}, \cite{WILK}).\\

In general we are concerned with groups $\mathcal{G}$ that act in an interesting way on relevant spaces (besides the natural action on $\mathcal{G}$ itself).  For instance, many of the above examples are {\it linear cocycles}, that is,   the elements of $\mathcal{G}$ act as linear maps of suitable vector bundles over the base space $\Omega$. A folklore approach allows to consider the elements of the general linear group $GL(d , \mathbb{R})$ as isometries of   the set of  positives matrices $Pos(d, \mathbb{R})$ endowed with a suitable metric (resulting into a non-positive curved space, see \cite{KAR_MAR}). With this in mind, we will be interested in cocycles taking values on the group $Isom(\mathcal{H})$ of isometries of a metric space $\mathcal{H}$ of negative  curvature (to be defined in a precise sense). In the context of a minimal base dynamics, in \cite{CONAPO} the authors show that a bounded continuous cocycle by isometries of a non-positively curved complete metric space   is always reducible, thus showing a general version of the Gottschalk-Hedlund Theorem. \\

In this work we will be placed in the Li\v vsic setting, that is, when $T$ is a hyperbolic homeomorphism (to be defined in a precise sense, see section \ref{CLT}). These maps have (among other features) a dense set of periodic orbits. Hence, it is interesting to consider the result of the cohomological equation on a periodic orbit. Let  $\omega$ be such that  there exists $n\in \mathbb{N}$ with  $T^n \omega=\omega$. By iterating the equation on such a point, we obtain a direct obstruction to the existence of  a (at least formal) solution of the cohomological equation

\begin{equation}\label{POO}
 \prod_{i=0}^{n-1} A(T^i\omega)
= \prod_{i=0}^{n-1} B(T^{i+1} \omega) \cdot B(T^i \omega)^{-1}
= B(T^n \omega)\cdot B(\omega)^{-1}=e_{_{\mathcal{G}}}.
\end{equation}
The {\it Liv\v sic problem} consists in determining whether the 
condition (\ref{POO}) is not only necessary
but also sufficient for $A$ being a coboundary. This terminology originates in the seminal work of
Liv\v sic \cite{LIV72}, who proved that this is the case whenever $\mathcal{G}$ is Abelian, $A$
is H\"older-continuous and $T$ is a topologically transitive hyperbolic diffeomorphism. Since then,
many extensions of this classical result have been proposed. Perhaps the most relevant is Kalinin's
recent version for $\mathcal{G} = \mathrm{GL} (d,\mathbb{C})$. 
\\

\noindent{\bf Theorem of Kalinin, see \cite{KAL10}.}
Let $T$ be a topologically transitive hyperbolic homeomorphism of a compact metric space $\Omega$. Let $A:\Omega\to GL(d, \mathbb{C})$ be an $\alpha$-H\"older function for which the 
 condition (\ref{POO}) holds. Then there exists an $\alpha$-H\"older function $B:\Omega\to GL(d, \mathbb{C})$ such that
for all $\omega \in \Omega$, $
A(\omega) = B(T\omega)\cdot  B(\omega)^{-1}.$\\

This impresive result is strongly based on the linear action of the matrices.  In the context of cocycles that take values on the group of diffeomorphisms of a manifold,  the answer 
to the Liv\v sic problem is
unclear. In the case of diffeomorphisms of a closed manifold (compact and without boundary), the recent result by Avila, Kocsard and Liu \cite{AVKOLI} is the most important step in the theory, extending the previous result by Kocsard and Potrie \cite{KOPO} on circle diffeomorphisms. \\

\noindent{\bf Theorem of Avila, Kocsard and Liu, see \cite{AVKOLI}. } Let $T$ be a topologically transitive hyperbolic homeomorphism of a compact metric space $\Omega$. Let $A:\Omega\to Diff^r(M)$ be an $\alpha$-H\"older cocycle, taking values on the group of diffeomorphisms of class $C^r$, $r>1$, of a smooth closed manifold $M$. If the 
 condition (\ref{POO}) holds then there exists a H\"older continuous map $B:\Omega\to Diff^r(M)$ such that $B(T\omega)\cdot B(\omega)^{-1}=A(\omega)$.  \\

 Our work  focuses on the study of the cohomological equation for cocycles of isometries of negatively curved metric spaces in the Li\v vsic setting.  The approach we propose is to take advantage of the fact that negatively curved  spaces, even though they are non compact, have a natural compactification as they are completed with a boundary at infinity. This boundary has a natural topology, that turns it into a compact space. For example, when we deal with a sectional-negatively curved riemannian manifold, the boundary at infinity is just a codimension one sphere (see \cite{EO}). In general, every isometry  extends to a homeomorphism of this boundary. In the context of cocycles by isometries, this produces a new cocycle by homeomorphisms of the boundary at infinity, for which the reducibility question can be  placed (compare with the Theorem by Avila, Kocsard and Liu above). The main goal of this work is to relate the reducibility of the cocycle by isometries with the reducibility of the cocycle by homeomorphisms induced on the boundary at infinity. \\
 
While the reducibility at the bounday at infinity is a more or less direct consequence of the reducibility of the original cocycle,  the other direction is more interesting and complicated. Although the notations and definitions of the hypotheses does not allow us to precisely state the main result in this introduction, we can paraphrase it as follows.\\
 
 \noindent{\bf  Result [see Theorem 5.3]. } Let $(T, A)$ be a continuous cocycle by isometries of a negatively curved complete metric space. Suppose that $T$ is hyperbolic, that $(T, A)$ verifies the condition (\ref{POO}) and that $(T, A)$ is H\"older. Let $(T,A^*)$ be the cocycle by homeomorphisms of the boundary at infinity. If $(T, A^*)$ is reducible by H\"older-Busseman functions then $(T, A)$ is reducible.\\

 \noindent{\bf Organization of the paper. } In section 2 we develop, in detail, the well-known fact that the existence of the solution to the cohomological equation is equivalent to the saturation by invariant sections of the phase space of a certain dynamical system (the skew-product). In section 3 we revisit the classical Li\v vsic theory in the real case (which serves as a context as well as being used in a key step in the proof of the main result of this work). In section 4 we review the geometric properties of the negative curvature spaces that we will consider (Gromov hyperbolic spaces). Section 5 is the core of this work, as it contains the precise statement of the main result, in addition to the description of the  H\"older conditions on the Busseman functions, which are part of the hypotheses of the main result . The section ends with the proof of the main result (Theorem 5.3). We complement the article with an appendix, in which we discuss the Hölder conditions on the Busseman functions of section 5 and their relationship with possible metrics at the boundary at infinity, in the case of negative curvature spaces in a more strict sense.

%%%%%%%%%%%%%%%%%%%%%%%%%%%%%%%%%%%%%%%%%%%
\section{Three equivalent problems}

This section is well known for the specialist. Nevertheless,   it worth some lines of paper to have a place with all the details of  these three equivalent problems related to reducibility of cocycles and it relation with skew-product dynamics.\\

Let's consider a continuous transformation $T:\Omega\to \Omega$, where $\Omega$ is a compact metric space, a complete metric space $\mathcal{H}$ and a function $A:\Omega\to Isom(\mathcal{H})$, that takes values on the space of isometries of $\mathcal{H}$. We write $d_{\mathcal{H}}(\cdot, \cdot), d_{\Omega}(\cdot, \cdot)$ for the distances in $\mathcal{H}$ and $\Omega$ respectively.\\

\noindent{\bf Definition 2.1.} We say that $A:\Omega\to Isom(\mathcal{H})$ is continuous for the topology of the uniform convergence on bounded sets when, given any bounded set $K\subset \mathcal{H}$, any $\omega_0\in \Omega$ and $\varepsilon >0$ there exists $\delta>0$ such that for every $h\in K$ we have
\[
d_{\Omega}(\omega_0,\omega )<\delta \Rightarrow d_{\mathcal{H}}(A(\omega)\cdot h, A(\omega_0)\cdot h)<\varepsilon.
\]

 Whenever $A$ is  continuous for the topology of the uniform convergence on bounded sets we will simply say that the pair $(T, A)$ is a {\it continuous cocycle by isometries} of the fiber $\mathcal{H}$ over the base space $\Omega$. \\
 
\noindent{\bf The cohomological equation.} The main problem we want to address in this work is about  the existence  of solutions   of the {\it cohomological equation}
\begin{displaymath}
(P1)\quad \left\{
\begin{array}{l}
    B(T\omega)\cdot B(\omega)^{-1}= A(\omega), \\
B:\Omega \to Isom(\mathcal{H})\ \textrm{is continuous for the pointwise topology.}
\end{array}
\right.
\end{displaymath}
We will also consider the following related problem 
\begin{displaymath}
(P1')\quad \left\{
\begin{array}{l}
    B(T\omega)\cdot B(\omega)^{-1}= A(\omega), \\
B:\Omega \to Homeo(\mathcal{H})\ \textrm{is continuous for the pointwise topology.}
\end{array}
\right.
\end{displaymath}
We write $Homeo(\mathcal{H})$ for the space of homeomorphisms of $\mathcal{H}$ with the metric topology, endowed with the pointwise topology.\\

\noindent{\bf Remark 2.2.} Given $G\in Isom(\mathcal{H})$ and a solution $B:\Omega\to Isom(\mathcal{H})$ to the problem $(P1)$, the function $B^{\star}=B\cdot G^{-1}$ also verifies $(P1)$. Hence, without loss of generality, given $\omega_0\in \Omega$ we can put $B(\omega_0)=id_{\mathcal{H}}$, the identity map of $\mathcal{H}$. An analogous remark also holds for the problem $(P1')$.\\

\noindent{\bf Lemma 2.3.} $Isom(\mathcal{H})$ is a closed subspace of $Homeo(\mathcal{H})$ for the pointwise topology.\\

\begin{proof} Let $(G_i)\subset Isom(\mathcal{H})$ be a  sequence that converges to $G\in Homeo(\mathcal{H})$. Pick any $h, h^*\in \mathcal{\mathcal{H}}$. The pointwise convergence implies that 
\[
\lim_{i\to \infty} d_{\mathcal{H}}(G_i(h), G_i(h^*))\longrightarrow d_{\mathcal{H}}(G(h), G(h^*)).
\]
Since $d_{\mathcal{H}}(G_i(h), G_i(h^*))=d_{\mathcal{H}}(h, h^*)$ we conclude that $G\in Isom(\mathcal{H})$.
\end{proof}

\noindent{\bf Lemma 2.4.} If $T$ admits a dense orbit then $(P1)\iff (P1')$.\\

\begin{proof} Suppose $B:\Omega\to Homeo(\mathcal{H})$ solves $(P1')$. Let $\omega_0\in \Omega$ whose $T-$orbit is dense. The relation $B(T\omega_0)\cdot B(\omega_0)^{-1}=A(\omega_0)$, its iterations and $B(\omega_0)=id_{\mathcal{H}}$ give raise to
\[
B(T^n\omega_0)=A(T^{n-1}\omega_0)\cdot A(T^{n-2}\omega_0)\cdots A(\omega_0).
\]
Hence  $B(T^n\omega_0)\in Isom(\mathcal{H})$ for every $n\in \mathbb{N}$. Given $\omega \in \Omega$, pick a sequence $(n_i)$ such that $T^{n_i}\omega_0\to \omega$. Continuity of $B$, together with the previous lemma allows to conclude that $B(\omega)=\lim_{i\to \infty}B(T^{n_i}\omega_0)$ belongs to $Isom(\mathcal{H})$.  The converse implication is direct.
\end{proof}

\noindent{\bf Skew product dynamical system.} Given a cocycle $(T, A)$ we can construct the following dynamics
\begin{eqnarray*}
F:\Omega\times \mathcal{H}&\longrightarrow& \Omega\times \mathcal{H}\\
(\omega, h)&\longmapsto& (T\omega, A(\omega)\cdot h),
\end{eqnarray*}
that is called the {\it induced skew-product}. The  simplest, notwithstanding key example, is the case $A(\omega)\equiv id_{\mathcal{H}}$. The skew-product $I_T(\omega, h)=(T\omega, h)$ is called {\it the fibred identity} over $T$.\\

Let $(T, A)$ be  a cocycle.  Given $G:\Omega \to Homeo(\mathcal{H})$  we construct
\begin{eqnarray*}
\mathcal{G}:\Omega\times \mathcal{H}&\longrightarrow& \Omega\times \mathcal{H}\\
(\omega, h)&\longmapsto& (\omega, G(\omega)\cdot h).
\end{eqnarray*}
We consider the following problem, written in the notation above:
\begin{displaymath}
(P2)\quad \left\{
\begin{array}{l}
    \mathcal{G}^{-1}\circ F\circ \mathcal{G}=I_{T}, \\
G:\Omega \to Homeo(\mathcal{H})\ \textrm{is continuous for the pointwise topology.}
\end{array}
\right.
\end{displaymath}
In that case we say that the skew product $F$ is conjugated to the fibred identity $I_{T}$ via a fiberwise preserving homeomorphism.\\

\noindent{\bf Lemma 2.5.} $(P1')\iff (P2)$.\\

\begin{proof} The second coordinate of the relation $\mathcal{G}^{-1}\circ F\circ \mathcal{G}=I_{T}$ is $G(T\omega)^{-1}\cdot A(\omega)\cdot G(\omega)=id_{\mathcal{H}}$. 
\end{proof}

\noindent{\bf Invariant sections. } In the notation above, the dynamics of the fibred identity $I_T$ is notting more than the consideration of many copies of the dynamics of $T$. In fact, 
$I^n_T(\omega, h)=(T^n\omega, h)$. In this way, the space $\Omega\times \mathcal{H}$ is decomposed into a disjoint union of $I_T$-invariant sections 
\[
\Omega\times \mathcal{H}=\bigcup_{h\in \mathcal{H}}\Omega\times \{h\}.
\]
Assuming the existence of a solution for $(P1')$ (and hence for $(P2)$), we will see that $F$ also induces a decomposition of $\Omega\times \mathcal{H}$ into a disjoint union of continuous and $F$-invariant sections. Let's define more precisely the elements of this claim. \\

Given $(\omega_0, h_0)\in \Omega\times \mathcal{H}$, we say that the function $s:\Omega\to \mathcal{H}$ is a {\it continuous invariant section passing through} $(\omega_0, h_0)$ if
\begin{enumerate}
    \item $s(\omega_0)=h_0$.
    \item For every $\omega\in \Omega$ we have
    \[
    A(\omega)\cdot s(\omega)=s(T\omega).
    \]
    
\end{enumerate}
\noindent{\it Example. } As commented before, if we have $A(\omega)\equiv Id_{\mathcal{H}}$ for every $\omega \in \Omega$ then $F=I_T$ and for every $(\omega_0, h_0)$ the constant section $s(\omega)=h_0$ is a continuous invariant section that passes through $(\omega_0, h_0)$.\\

\noindent{\it Example. } Let $(T, A)$ a continuous cocycle by isometries and such that $(P1')$ has a continuous solution $B:\Omega\to Homeo(\mathcal{H})$  (or equivalently $(P2)$ has a solution). We know that $\mathcal{G}(\omega, h)=(T\omega, B(\omega)\cdot h)$ conjugates $F$ with $I_T$. Hence, for $(\omega_0, h_0)$ we define 
\[
s_{\omega_0, h_0}(\omega)=B(\omega)\cdot B(\omega_0)^{-1}\cdot h_0.
\]
Continuity is direct from pointwise continuity of $B$. We can verify the  conditions:
\begin{enumerate}
    \item $s_{\omega_0,h_0}(\omega_0)=B(\omega_0)\cdot B(\omega_0)^{-1}\cdot h_0=h_0.$
    \item \begin{eqnarray*}
    A(\omega)\cdot s_{\omega_0,h_0}(\omega)&=&A(\omega)\cdot B(\omega)\cdot B(\omega_0)^{-1}\cdot h_0\\
    &=&B(T\omega)\cdot B(\omega_0)^{-1}\cdot h_0\\
    &=&s_{\omega_0, h_0}(T\omega). 
    \end{eqnarray*}
\end{enumerate}
This example, and it reciprocal, is the key of this work, since it relates the solution of the cohomological equation  with the existence of continuous invariant sections for the induced skew product, that passes through every point of the product space. Indeed, we will consider the following problem.\\

Given a continuous cocycle by isometries $(T, A)$,  and the corresponding skew product dynamics $F(\omega, h)=(T\omega, A(\omega)\cdot h)$, we ask whether or not, the following holds
\begin{displaymath}
(P3)\quad \left\{
\begin{array}{l}
    \textrm{For every $(\omega_0, h_0)$ there exists a continuous invariant section $s_{\omega_0, h_0}$}\\ 
    \textrm{that passes through $(\omega_0, h_0)$}.\\
     \textrm{For every $\omega_0, \omega\in \Omega$ the map $ h_0\mapsto s_{\omega_0, h_0}(\omega)$ is continuous.}
\end{array}
\right.
\end{displaymath}
\noindent{\bf Lemma 2.6. } Let $(T, A)$ be a continuous cocycle by isometries. Let $(\omega_0, h_0)\in \Omega\times \mathcal{H}$. If $s, \tilde s$ are two continuous invariant sections that passes through $(\omega_0, h_0)$ and $T$ admits a dense orbit then $s=\tilde s$.\\

\begin{proof} Assume there exists $\omega_1\in \Omega$ such that $\tilde s(\omega_1)\neq s(\omega_1)$. Let $\omega_{*}$ with dense $T$-orbit. Taking $T^{n_*}\omega_*$ close enough to $\omega_1$, we can find a positive integer $n_*\in \mathbb{N}$ such that
\begin{equation}\label{eq1}
d_{\mathcal{H}}(\tilde s(T^{n_*}\omega_*), s(T^{n_*}\omega_*))>\frac{1}{2}d_{\mathcal{H}}(\tilde s(\omega_1), s(\omega_1))>0.
\end{equation}
Let $(n_j)$ be a sequence of positive integers such that $T^{n_*+n_j}\omega_*\to \omega_0$. Since $s, \tilde s$ are continuous and both passe through $(\omega_0, h_0)$, we verify that
\begin{eqnarray*}
\lim_{j\to \infty} s(T^{n_*+n_j}\omega_*)&=&s(\omega_0)=h_0, \\
\lim_{j\to \infty} \tilde s(T^{n_*+n_j}\omega_*)&=&s(\omega_0)=h_0.
\end{eqnarray*}
Then we obtain
\begin{equation}\label{eq2}
    \lim_{j\to \infty} d_{\mathcal{H}}(\tilde s(T^{n_*+n_j}\omega_*), s(T^{n_*+n_j}\omega_*))=0.
\end{equation}
Since both $s, \tilde s$ are invariant by $F$ we have
\begin{eqnarray*}
s(T^{n_*+n_j}\omega_*)&=&A^{n_j}(T^{n_*}\omega_*)\cdot s(T^{n*}\omega_*), \\
\tilde s(T^{n_*+n_j}\omega_*)&=&A^{n_j}(T^{n_*}\omega_*)\cdot s(T^{n*}\omega_*).
\end{eqnarray*}
As $A^{n_j}(T^{n_*}\omega_*)\in Isom(\mathcal{H})$ for every $j$, we obtain
\[
d_{\mathcal{H}}(\tilde s(T^{n_*+n_j}\omega_*), s(T^{n_*+n_j}\omega_*))=d_{\mathcal{H}}(\tilde s(T^{n_*}\omega_*), s(T^{n_*}\omega_*)),
\] 
which is incompatible with (\ref{eq1}) and (\ref{eq2}), concluding the proof.
\end{proof}

\noindent{\bf Remark 2.7.} Note that the problem $(P3)$ can be also stated for cocycles of homeomorphisms $A:\Omega \to \mathcal{G}$, where $\mathcal{G}$ is a topological space. The formal construction of the skew-product is  the same. \\

\noindent{\bf Proposition 2.8. } Let $(T, A)$ be a continuous cocycle by isometries. If  $T$ admits a dense orbit then $(P2)\iff (P3)$.\\

\begin{proof} We need to show that $(P3)$ implies $(P2)$, since the other direction was already discussed in the examples. Fix $\omega_0\in \Omega$. For every $\omega\in \Omega$ we need to define $B(\omega)\in Homeo(\mathcal{H})$ solving $(P2)$. We define
\[
B(\omega)\cdot h= s_{\omega_0, h} (\omega).
\]
Let's verify the conditions on $(P2)$.  First, we notice that $B(\omega_0)\cdot h=s_{\omega_0, h}(\omega_0)=h$ and hence $B(\omega_0)=id_{\mathcal{H}}$. We also have
\begin{eqnarray*}
B(T\omega)\cdot h&=& s_{\omega_0, h}(T\omega)\\
&=& A(\omega)\cdot s_{\omega_0, h}(\omega)\\
&=&A(\omega)\cdot B(\omega)\cdot h,
\end{eqnarray*}
thus, $B(T\omega)=A(\omega)\cdot B(\omega)$. We still have to verify that $B(\omega)\in Homeo (\mathcal{H})$. The function $h\mapsto s_{\omega_0, h}(\omega)=B(\omega)\cdot h$ is continuous by $(P3)$. Thanks to the uniqueness of the invariant section, the inverse of $B(\omega)$ can be easily computed as
\[
B(\omega)^{-1}\cdot g= s_{\omega, g}(\omega_0),
\]
which is also continuous. The pointwise continuity of $\omega\mapsto B(\omega)\in Homeo(\mathcal{H})$ comes from the continuity of $\omega \mapsto s_{\omega_0, h}(\omega)$.
\end{proof}

In the previous proof, the choice of $\omega_0$ is arbitrary, giving account of the non-uniqueness of the solution to the cohomological equation. We summarize the discussion of the current section in the following theorem.\\

\noindent{\bf  Theorem 2.9 [folklore].} Let $(T,A)$ be a continuous cocycle by isometries and such that $T$ admits a dense orbit. The following statements are equivalent.\\

\noindent{\bf (P1)} The cohomological equation
\[
B(T\omega)\cdot B(\omega)^{-1}=A(\omega)
\]
admits a  solution $B:\Omega\to Isom(\mathcal{H})$, that is continuous for the pointwise topology.\\

\noindent{\bf (P2)} There exists $G:\Omega\to Homeo(\mathcal{H})$ continuous for the pointwise topology, and such that
\begin{eqnarray*}
\mathcal{G}:\Omega\times \mathcal{H}&\to&\Omega\times\mathcal{H}\\
(\omega, h)&\mapsto& (\omega, G(\omega)\cdot h)
\end{eqnarray*}
verifies
\[
\mathcal{G}^{-1}\circ F\circ \mathcal{G}=I_T, 
\]
where $F(\omega, h)=(T\omega, A(\omega)\cdot h)$ is the  skew product induced by $(T, A)$ and $I_T(\omega, h)=(T\omega, h)$ is the fibred identity. \\

\noindent{\bf (P3)} For every $(\omega_0, h_0)\in \Omega\times \mathcal{H}$ there exists a continuous section $s:\Omega \to \mathcal{H}$ that verifies
\begin{enumerate}
    \item $s(\omega_0)=h_0$, 
    \item $A(\omega)\cdot s(\omega)=s(T\omega)$ for every $\omega\in \Omega$. 
\end{enumerate}
Moreover the map $ h_0\mapsto s$ is continuous for the pointwise topology.
\section{Classic Li\v vsic Theorem}\label{CLT}
In this section we revisit the Classic Li\v vsic Theorem on cohomological equations (see \cite{LIV72}). Let $\Omega$ be a complete metric space. We say that a homeomorphism $T:\Omega \to \Omega$ is {\it hyperbolic} if the following two conditions holds.
\begin{enumerate}
\item $T$ admits a dense orbit.
    \item Let   $x, y$ be any pair of points of $\Omega$. We 
say that the orbit segments  $x, Tx, \dots, T^kx$  and  
$y, Ty, \dots, T^ky$  are {\it exponentially $\delta$-close 
with exponent  $\lambda>0$}  if for every $j=0, \dots, k$,
\[
d_{\Omega}(T^jx,T^jy)\leq \delta e^{-\lambda\min\{j, k-j\}}.
\]
We say that $T$ satisfies the {\it closing property} if there exist 
$c, \lambda, \delta_0>0$ such that for every $x\in \Omega$ and $k\in \mathbb{N}$ 
so that $d_{\Omega}(x, T^kx)<\delta_0$, there exists a point $p\in X$ with 
$T^kp=p$ so that letting $\delta := c  d_{\Omega}(x, T^kx)$, the orbit segments 
$x, Tx, \dots, T^kx$ and $p, Tp, \dots, T^kp$ are exponentially $\delta$-close 
with exponent $\lambda$ and there exists a point $y \in \Omega$ such that for every $j=0,\dots, k$,
\[
d_{\Omega}(T^jp, T^jy)\leq \delta e^{-\lambda j} \quad \textrm{and} 
\quad d_{\Omega}(T^jy, T^jx)\leq \delta e^{-\lambda(n-j)}.
\]
    
\end{enumerate}

Important examples of maps satisfying these properties are hyperbolic diffeomorphisms of compact manifolds.\\

Given a measurable function $\psi:\Omega\to \mathbb{R}$, we consider the real valued cohomological equation 
\begin{equation}\label{Livsic_real}
u(T\omega)-u(\omega)=\psi(\omega).
\end{equation}
By taking integrals at both sides of this equation, we discover that an immediate obstruction in order to obtain an, at least, measurable solution $u:\Omega\to \mathbb{R}$ is of that for every $T$-invariant measure $\mu$ we should have $\int_{\Omega}\psi(\omega)d\mu=0$. In the particular case of a periodic point this integral condition reads
\begin{displaymath}
(PPO)_{\mathbb{R}}\quad \left\{
\begin{array}{c}
    \omega^*\in \Omega\ \textrm{and}\ n^*\in \mathbb{N} \ \textrm{such that}\ T^{n^*}\omega^*=\omega^* \\
    \Downarrow\\
    \sum_{j=0}^{n^*-1}\psi(T^j\omega^*)=0.
\end{array}
\right.
\end{displaymath}
\noindent{\bf Theorem 3.1 see [Li\v vsic \cite{LIV72}].} Let $T:\Omega\to \Omega$ be a hyperbolic dynamical system defined on a complete metric space $\Omega$. Given a H\"older function $\psi:\Omega\to \mathbb{R}$ verifying the $(PPO)_{\mathbb{R}}$ condition, there exists a H\"older solution $u:\Omega \to \mathbb{R}$ to the cohomologival equation (\ref{Livsic_real}).\\

\noindent{\bf Remark 3.2.} Given two continuous solutions $u_1, u_2$ to (\ref{Livsic_real}), the difference $f=u_1-u_2$ verifies $f(T\omega)=f(\omega)$ for every $\omega\in \Omega$. Since $T$ admits a dense orbit we conclude that $f$ is constant. In other words, any two continuous solutions to (\ref{Livsic_real}) differ by a constant. 

\section{Gromov Hyperbolic metric spaces}
This section collects some geometric elements that we use in the setting of our main theorem. The main reference is the almost comprehensive book by Bridson and H\"afliger \cite{BRIDSON}, chapter III. We develop some of the lemmas in order to bring attention to some relevant ideas and techniques. \\

Given a proper metric space $\mathcal{H}$, with distance denoted by $d_{\mathcal{H}}$, recall that a {\it geodesic} is an isometric map $\gamma:[a,b]\subset \mathbb{R} \to \mathcal{H}$. The space $\mathcal{H}$ is called geodesic if every pair of points can be joined by a geodesic; if this geodesic is unique we call $\mathcal{H}$ a {\it unique geodesic space}. A geodesic triangle (or simply a triangle) is a set consisting of the union of three geodesics that join three points $x,y,z$ between each other. Points are called vertices and the geodesic that join two vertices, $x,y$, denoted by $[x,y]$, is called a side of the triangle.\\

\noindent{\bf Definition 4.1.} Let $\delta>0$. A geodesic triangle  $\triangle ABC$ is said to be $\delta$-\textit{slim} if each of its sides is contained in a $\delta$ neighborhood of the union of the other two sides. That is to say, for every point $p\in [A,B]$ there exists $q\in [B,C]\cup [A,C]$ such that $d_{\mathcal{H}}(p,q)\leq \delta$.  A  unique geodesic space $\mathcal{H}$ is said to be {\it Gromov hyperbolic} if there exists $\delta>0$ so that every triangle is $\delta$-slim.\\

A first  example of these spaces are real trees (non directed graphs in which two vertices are joined by exactly one path). Classic hyperbolic space $\mathbb{H}^n$ is also Gromov hyperbolic, and, as a consequence, $\textsc{CAT}(\kappa)$ spaces are Gromov hyperbolic, for $\kappa<0$.  A lot of interesting examples of Gromov hyperbolic spaces were given by Gromov in \cite{Gromov1987}. \\

\noindent{\bf Definition 4.2.}
 We will say two geodesic rays $\gamma_1,\gamma_2:[0,\infty]\to \mathcal{H}$,  are equivalent if there exists $K>0$ such that  $\sup_{t\geq 0} d_{\mathcal{H}}(\gamma_1(t),\gamma_2(t))\leq K$. The \emph{Gromov Boundary} $\partial \mathcal{H}$ is defined as the equivalence classes among all geodesic rays of the space. Given a geodesic $\gamma:[-\infty, \infty]\to \mathcal{H}$ we write $\gamma(+\infty)=[\gamma]$ and $\gamma(-\infty)=[t\mapsto \gamma(-t)]$, and we call it the {\it ends of } $\gamma$.\\
 
\noindent{\bf Proposition 4.3 [see \cite{BRIDSON} III.3.2.]. } Let $\mathcal{H}$ be a (proper unique geodesic) Gromov hyperbolic  space. The boundary $\partial \mathcal{H}$ is {\it visible}, that is, given $x\in \mathcal{H}$ and a boundary point  $\alpha \in \partial \mathcal{H}$, there exists a unique  geodesic ray $\gamma:[0, \infty]\to \mathcal{H}$ such that $\gamma(+\infty)=\alpha$ and $\gamma(0)=x$. Moreover, given two boundary points $\alpha, \beta \in \partial{\mathcal{H}}$ there exists a geodesic $\gamma:[-\infty, \infty]\to \mathcal{H}$ such that $\gamma(+\infty)=\alpha$ and $\gamma(-\infty)=\beta$.  \\

\noindent{\bf Definition 4.4. } A unique geodesic Gromov hyperbolic space will be called {\it a unique visibility Gromov space} whenever the geodesic joining every pair of boundary points is unique. \\

A {\it generalized geodesic ray} is a geodesic ray $\gamma:[0, \infty]\to \mathcal{H}$ or a path that is a geodesic until a certain point and then constant. The next definition provides a topology for $\partial \mathcal{H}$. From now on we will assume that $\mathcal{H}$ is a proper unique visibility Gromov hyperbolic space.\\

\noindent{\bf Definition 4.5.}
Let $p\in \mathcal{H}$. We say that a sequence $(x_n)\subset \mathcal{H}$ converges to a point $x\in \mathcal{H}\cup \partial \mathcal{H}$, if there exist generalized geodesic rays $(\gamma_n)$ such that $\gamma_n(0)=p$ and $\gamma_n(\infty)=x_n$ and verifying that every subsequence of $(\gamma_n)$ contains a subsubsequence that converges, in the compact-uniform topology, to a generalized ray $\gamma$ with $\gamma(+\infty)=x$.\\

Given an isometry $A\in Isom(\mathcal{H})$ and a geodesic ray $\gamma:[0, \infty]\to \mathcal{H}$, the image $A\cdot \gamma$ is a geodesic ray and hence the action of $A$ on $\mathcal{H}$ can be extended to an action $A^*:\partial \mathcal{H}\to \partial \mathcal{H}$ by considering the equivalence classes $A^*\cdot [\gamma]=[A\cdot \gamma]$. Since $A$ is an isometry, $A^*$ is well defined.\\

\noindent{\bf Lemma 4.6. } Given two isometries $A, B$ we have $(A\circ B)^*=A^*\circ B^*$.\\

\noindent{\bf Lemma 4.7.}  $A^*=id_{\partial \mathcal{H}}$ if and only if  $A=id_{\mathcal{H}}$.

\begin{proof} Let $x\in \mathcal{H}$. Choose $\alpha_1\in \partial \mathcal{H}$ and $\gamma_1$ the geodesic such that $\gamma_1(0)=x$ and $\gamma_1(+\infty)=\alpha_1$. Let's call $\beta_1=\gamma_1(-\infty)$ and pick $\alpha_2\in \partial \mathcal{H}$ different to $\alpha_1$ or $\beta_1$. Define $\gamma_2$ as the unique geodesic such that $\gamma_2(0)=x$ and $\gamma_2(+\infty)=\alpha_2$. Notice that $\gamma_1\neq \gamma_2$ and hence $x$ is the unique point in $\gamma_1\cap \gamma_2$. Since $A^*\cdot \alpha_1=\alpha_1$, $A^*\cdot \beta_1=\beta_1$, and $A^*\cdot \gamma_1$ is a geodesic, we conclude that $A^*\cdot \gamma_1=\gamma_1$ and analogously $A^*\cdot \gamma_2=\gamma_2$. As $A^*\cdot x$ equals $A^*\cdot \gamma_1\cap A^*\cdot \gamma_2$ we conclude $A^*\cdot x=x$.
\end{proof}

\noindent{\bf Lemma 4.8. }  If $A\in Isom(\mathcal{H})$ then $A^*:\partial\mathcal{H}\to \partial\mathcal{H}$ is a homeomorphism. \\

%\noindent{\bf Theorem [Cartan-Hadamard]  G1.} Let $\mathcal{H}$ be a Hadamard complete metric space. For every pair $(x, y)\in \mathcal{H}$ there exists a unique geodesic segment $\gamma_{x, y}:[0, d_{\mathcal{H}}(x, y)]\to \mathcal{H}$ such that $\gamma_{x, y}(0)=x, \gamma_{x, y}(d_{\mathcal{H}}(x, y))=y$. This geodesic segment minimizes the distance between any pair of its points. Moreover, the map $(x, y)\mapsto \gamma_{x, y}$ is continuous.\\

\noindent{\bf Lemma 4.9.} Let $\alpha_0, \beta_0 \in \partial \mathcal{H}$ and  the geodesic  $\gamma_0:[-\infty, \infty]\to \mathcal{H}$ such that $\gamma_0(-\infty)=\beta_0\in \mathcal{H}$ and $\gamma_0(+\infty)=\alpha_0$. Consider two sequences $(\alpha_n, \beta_n)_{n\in \mathbb{N}}$ in $\partial \mathcal{H}$ such that $\beta_n\to \beta$ and $\alpha_n\to \alpha$. Define $\gamma_n$ as the complete geodesic given by Proposition  4.3., such that $\gamma_n(-\infty)=\beta_n$ and $\gamma_n(+\infty)=\alpha_n$. For any $r>0$, $s\in \mathbb{R}$ there exists $N\in \mathbb{N}$ such that for  $n>N$ 
\[
\gamma_n\cap B(\gamma_0(s), r)\neq \emptyset.
\]

\noindent{\bf Lemma 4.10.} Let $h\in \mathcal{H}$ and $r>0$ small enough. Then for every pair of points $x, y\in B(h, r)$ the geodesic segment $\gamma_{x, y}\subset B(h, r)$. That is, the ball $B(h, r)$ is geodesically convex.\\

\noindent{\bf Lemma 4.11.} Let be $\alpha\in \partial\mathcal{H}$. For every $h\in \mathcal{H}$ let $\gamma_h:[0, \infty]\to \mathcal{H}$ be the unique geodesic joining $h$ with $\alpha$. Extend  $\gamma_h$ to a complete geodesic and denote by $\beta_h$ the point $\gamma_h(-\infty)\in \partial \mathcal{H}$. Then the map $h\mapsto \beta_h$ is continuous.

\subsection{Busseman functions}
Let $\mathcal{H}$ be a unique geodesic Gromov hyperbolic space. Given $\alpha\in \partial \mathcal{H}$ and a point $p\in\mathcal{H}$, we define the Busseman function $b_{p, \alpha}:\mathcal{H}\to \mathbb{R}$, {\it in the direction} $\alpha$  and {\it with base point} $p$ as
\[
b_{p, \alpha}(h)=\lim_{n\to \infty} d_{\mathcal{H}}(x_n, h)-d_{\mathcal{H}}(x_n, p),
\]
where $(x_n)\subset \mathcal{H}$ is any sequence such that $x_n\to \alpha$. The convergence and independence on $(x_n)$ relies on the triangle inequality and the hyperbolicity of $\mathcal{H}$.\\

\noindent{\bf Definition 4.12.} Given any real number $r\in \mathbb{R}$, the level set $b_{p, \alpha}^{-1}(r)$ is called an {\it horosphere} centered at $\alpha$.\\

\noindent{\bf Lemma 4.13.} Let $\gamma:\mathbb{R}\to \mathcal{H}$ be a geodesic. If  $\gamma(+\infty)=\alpha$ and $b_{p, \alpha}(\gamma(0))=0$ then $b_{p, \alpha}(\gamma(s))=-s$ for every $s\in \mathbb{R}$.

\begin{proof} Note that $x_n=\gamma(n)$ verifies $x_n\to \alpha$. We have
\begin{eqnarray*}
d_{\mathcal{H}}(\gamma(n), \gamma(s))-d_{\mathcal{H}}(\gamma(n), p)&=&n-s-\left(d_{\mathcal{H}}(\gamma(n), \gamma(0))+d_{\mathcal{H}}(\gamma(n), p)-d_{\mathcal{H}}(\gamma(n), \gamma(0))\right)\\
&=&-s+\left(d_{\mathcal{H}}(\gamma(n), \gamma(0))-d_{\mathcal{H}}(\gamma(n), p)\right).
\end{eqnarray*}
Taking $n\to \infty$ above we obtain
\begin{eqnarray*}
b_{p, \alpha}(\gamma(s))&=&-s-b_{p, \alpha}(\gamma(0))\\
&=&-s.
\end{eqnarray*}
In the same lines, we can deduce $b_{p, \alpha}(\gamma(s))=b_{p, \alpha}(\gamma(r))+r-s$ for every $r, s\in \mathbb{R}$.
\end{proof}

\noindent{\bf Lemma 4.14. } Let $A\in Isom(\mathcal{H})$, $\alpha \in \partial \mathcal{H}$ and $p\in \mathcal{H}$. Then the following holds
\[
b_{p, A^*\cdot \alpha}(A\cdot p)=- b_{p, \alpha}(A^{-1}\cdot p)=b_{A^{-1}\cdot p, \alpha}(p).
\]
\begin{proof} Let $(x_n)\subset \mathcal{H}$ such that $x_n\to \alpha$. Hence $A\cdot x_n\to A^*\cdot \alpha$. Then
\begin{eqnarray*}
b_{p, A^*\cdot \alpha}(A\cdot p)&=&\lim_{n\to \infty} d_{\mathcal{H}}(A\cdot x_n, A\cdot p)-d_{\mathcal{H}}(A\cdot x_n, p)\\
&=&\lim_{n\to \infty}-\left(d_{\mathcal{H}}(x_n, A^{-1}\cdot p)-d_{\mathcal{H}}( x_n, p)\right)\\
&=&-b_{p, \alpha}(A^{-1}\cdot p).
\end{eqnarray*}
\end{proof}
\noindent{\bf Lemma 4.15. } If $\alpha \in \partial \mathcal{H}$ and $p\in \mathcal{H}$ then for every $h, g \in \mathcal{H}$ one has
\[
|b_{p, \alpha}(h)-b_{p, \alpha}(g)|\leq d_{\mathcal{H}}(h, g).
\]
\begin{proof} Given $x_n\to \alpha$, one knows that

\begin{eqnarray*}
|d_{\mathcal{H}}(x_n, h)-d_{\mathcal{H}}(x_n, p)-d_{\mathcal{H}}(x_n, g)+d_{\mathcal{H}}(x_n,p)|&=&|d_{\mathcal{H}}(x_n, h)-d_{\mathcal{H}}(x_n, g)|\\
&\leq& d_{\mathcal{H}}(h, g).
\end{eqnarray*}
Taking $\lim_{n\to \infty}$ we complete the claim.
\end{proof}

\noindent{\bf Lemma 4.16.} Given $p_1, p_2, h\in \mathcal{H}$ and $\alpha\in \partial \mathcal{H}$ we have
\[
b_{p_1, \alpha}(h)=b_{p_2, \alpha}(h)+b_{p_1, \alpha}(p_2).
\]
\begin{proof} Given $x_n\to \alpha$ one knows
\[
d_{\mathcal{H}}(x_n, h)-d_{\mathcal{H}}(x_n, p_1)=d_{\mathcal{H}}(x_n, h)-d_{\mathcal{H}}(x_n, p_2)+d_{\mathcal{H}}(x_n, p_2)-d_{\mathcal{H}}(x_n, p_1), 
\]
which, taking the limit $n\to \infty$ yields the result.
\end{proof}
\section{Boundary solutions imply interior solutions}

Given a cocycle by isometries $(T, A)$ we can easily induce a cocycle by homeomorphisms of the boundary $\partial \mathcal{H}$,  $A^*:\Omega\to Homeo(\partial \mathcal{H})$ defined as $A^*(\omega)=(A(\omega))^*$. We call $(T, A^*)$ the {\it boundary cocycle} induced by $(T, A)$. The following is direct from the definition of the topology on $\partial\mathcal{H}$.\\

\noindent{\bf Lemma 5.1.} If $(T,A)$ is  continuous then $(T,A^*)$  is continuous for the pointwise topology.\\ 

This lemma allows to fit the boundary cocycle $(T, A^*)$ into the framework of Problem (P1'), that is, to look for a pointwise continuous solution $B:\Omega\to Homeo(\partial \mathcal{H})$ to the cohomological equation $\overline{B}(T\omega)\cdot \overline{B}(\omega)^{-1}=A^*(\omega)$. \\

\noindent{\bf Remark 5.2. } This work is entirely devoted to relate the solutions to the cohomological equations for the cocycle of isometries $(T, A)$ and the boundary cocycle $(T, A^*)$. The  direction from interior solutions the problem $(P1) (or  (P2)$ or $(P3)$) to a solution of the corresponding problem for the boundary cocycle can be addressed in different ways. For example,  if $(T, A)$ is a continuous cocycle by isometries such that there exists a pointwise continuous solution $B:\Omega\to Isom(\mathcal{H})$ to the cohomological equation $B(T\omega)\cdot B(\omega)^{-1}=A(\omega)$, then Lemma 4.13. gives 
 $
 B^*(T\omega)\cdot B^*(\omega)^{-1}=A^*(\omega).
 $
 This relates $(P1)$ for $(T, A)$ and $(P1')$ for $(T, A^*)$. Nevertheless, pointwise continuity of $B^*$ can not be ensured in all cases. The existence of solutions to the problem $(P3)$ for the cocycle $(T, A)$ gives a solution to $(P3)$ for the cocycle $(T, A^*)$ in the following way. Assume $(T, A)$ has solutions to $(P3)$, that is, through every point $(\omega_0, h_0)$ passes a continuous $A$-invariant section $s_{\omega_0, h_0}(\omega)$, and such that for fixed $(\omega_0, \omega) $ the map $h\mapsto s_{\omega_0, h}(\omega)$ is continuous. We want to show that $(T, A^*)$ also has solutions for the corresponding $(P3)$. Fix $p\in \mathcal{H}$ and any pair $(\omega_0, \alpha_0)\in \Omega\times \partial \mathcal{H}$. Let $\gamma_0:\mathbb{R}^+\to \mathcal{H}$ be the geodesic ray that connects $p$ to $\alpha_0$ and let $\overline{p}=\gamma_0(1)$. Let $\gamma_{\omega}:\mathbb{R}^+\to \mathcal{H}$ be the geodesic ray that starts at $s_{\omega_0, p}(\omega)$ and passes through $s_{\omega_0, \overline{p}}(\omega)$. Define $s^*_{\omega_0, \alpha_0}(\omega)=\gamma_{\omega}(+\infty)$. It is easy to show that $s^*_{\omega_0, \alpha_0}$  is a continuous  $A^*$-invariant section that passes through $(\omega_0, \alpha_0)$. \\
 
The next is the main result in this work,  showing that, under the setting of the Li\v vsic  result for cohomological equations, somehow the converse of this relation holds. Definitions of the  required  {\it \`a la} H\"older conditions and suitable (PPO) conditions are described in the next sections \ref{52}, \ref{CONHOL}.   \\

\noindent{\bf Theorem 5.3. } Let $(T, A)$ be a cocycle by isometries of  a unique visibility Gromov hyperbolic  space $\mathcal{H}$. If $(T, A)$ satisfies 
\begin{itemize}
    \item $T$ is hyperbolic (in the sense of section  \ref{CLT}),
    \item $(T, A)$ is $\tau$-H\"older for some $\tau>0$,
    \item $(T, A)$ verifies the $(PPO)$ condition,
    \item $(T, A^*)$ verifies $(P3)$ with $\tau$-H\"older-Busemann invariant sections,
\end{itemize}
then $(T, A)$ verifies $(P3)$ (and equivalently $(P2)$ and $(P1)$).\\

\noindent{\bf Remark 5.4.} Lemmas 4.6., 4.7. imply that $(T, A)$ verifies the $(PPO)$ condition if and only if $(T, A^*)$ verifies the $(PPO)$ condition.

\subsection{Periodic Point Obstruction}\label{52}
By iterating the cohomological equation, we see that a necessary condition on $(T,A)$ in order to admit a solution to the cohomological equation $(P1')$ is the following {\it Periodic Point Obstruction} 
\begin{displaymath}
(PPO)\quad \left\{
\begin{array}{c}
    \omega^*\in \Omega\ \textrm{and}\ n^*\in \mathbb{N} \ \textrm{such that}\ T^{n^*}\omega^*=\omega^* \\
    \Downarrow\\
    A^n(\omega^*)=id_{\mathcal{H}}
\end{array}
\right.
\end{displaymath}
Notice that this condition generalizes the condition $(PPO)_{\mathbb{R}}$ in section \ref{CLT}. Indeed, the translations $G(\omega)\cdot t=t+\psi(\omega)$ are isometries of the real line $\mathbb{R}$ endowed with the usual $|\cdot - \cdot|$ distance, and  iterates $G^n(\omega)\cdot t$ reduce to the translation by the Birkhoff sum $\sum_{j=0}^{n-1}\psi(T^{j}\omega)$. Hence for a periodic point $T^{n*}(\omega^*)=\omega^*$ we have
\[
G^{n^*}(\omega^*)=id_{\mathbb{R}}\iff \sum_{j=0}^{n^*-1}\psi(T^j\omega^*)=0.
\]

The PPO condition on $(T,A)$ can be rewritten in the following way in terms of the skew product $F$. 
\begin{displaymath}
(PPO')\quad \left\{
\begin{array}{c}
    \omega^*\in \Omega\ \textrm{and}\ n^*\in \mathbb{N} \ \textrm{such that}\ T^{n^*}\omega^*=\omega^* \\
    \Downarrow\\
    F^{n^*}(\omega^*, h)=(\omega^*, h)\ \forall h\in \mathcal{H}
\end{array}
\right.
\end{displaymath}
\noindent{\bf Lemma 5.5.} $(PPO)\iff (PPO')$.\\

\noindent{\bf Lemma 5.6. } If $(T,A)$ verifies the $(PPO)$ condition then the induced cocycle on the boundary $(T, A^*)$ also verifies the $(PPO)$ condition.\\

\begin{proof} It is clear from $(A^n)^{*}=(A^{*})^n$.
\end{proof}
\subsection{H\"older conditions}\label{CONHOL} In this section we introduce the suitable H\"older like conditions that we use in this work.\\

\noindent{\bf Definition 5.7.} We say that  a cocycle by isometries $(T, A)$ is $\tau$-H\"older, $\tau>0$, if for every bounded set $K\subset \mathcal{H}$ there exists $C_K>0$ such that for every $\omega_1, \omega_2\in \Omega$, $h\in K$ one has
\[
d_{\mathcal{H}}(A(\omega_1)\cdot h, A(\omega_2)\cdot h)\leq C_Kd_{\Omega}(\omega_1, \omega_2)^{\tau}.
\]
\noindent{\bf Lemma 5.8. } Let $(T, A)$ be a $\tau$-H\"older cocycle by isometries and  $p\in \mathcal{H}$. There exists a constant $C_p>0$ such that for every  $\omega_1, \omega_2\in \Omega$
\[
d_{\mathcal{H}}(A(\omega_1)^{-1}\cdot p, A(\omega_2)^{-1}\cdot p)\leq C_pd_{\Omega}(\omega_1, \omega_2)^{\tau}.
\]
\begin{proof} Since $A(\omega_1)$ is an isometry, for every $h$ we have
\[d_{\mathcal{H}}(A(\omega_1)\cdot h, A(\omega_2)\cdot h)=d_{\mathcal{H}}( h, A(\omega_1)^{-1}\cdot A(\omega_2)\cdot h)
\]
Notice that the set $K_p=\{A^{-1}(\omega)\cdot p\ | \ \omega \in \Omega\}$ is bounded (since $d_{\mathcal{H}}(A(\omega)^{-1}\cdot p, p)=d_{\mathcal{H}}(p, A(\omega)\cdot p)$), and hence we can take a uniform constant $C_p=C_{K_p}$ such that for every $h=A(\omega_2)^{-1}\cdot p\in K_p$ we have the desired inequality.
\end{proof}

\noindent{\bf Definition 5.9.} We say that a section $\alpha:\Omega\to \partial \mathcal{H}$ is $\tau$-H\"older$-$Busemann, $\tau>0$, if for every bounded set $K\subset \mathcal{H}$ there exists $D_K>0$ such that for every $\omega_1, \omega_2\in \Omega$, $h\in K$ one has
\[
\left| b_{p, \alpha(\omega_1)}(h)-b_{p, \alpha(\omega_2)}(h)\right|\leq D_Kd_{\Omega}(\omega_1, \omega_2)^{\tau}.
\]
\noindent{\bf Lemma 5.10. } Let $\alpha:\Omega\to \partial \mathcal{H} $ be a $\tau$-H\"older$-$Busemann section  and $p\in \mathcal{H}$. There exists a constant $D_p>0$ such that for every  $\omega_1, \omega_2\in \Omega$
\[
\left| b_{p, \alpha(\omega_1)}(A(\omega_1)^{-1}\cdot p)-b_{p, \alpha(\omega_2)}(A(\omega_1)^{-1}\cdot p)\right|\leq D_pd_{\Omega}(\omega_1, \omega_2)^{\tau}.
\]
\subsection{The induced cocycle in the space of horospheres}
Assume that there exists a section $\alpha:\Omega\to \partial \mathcal{H}$ that is continuous and invariant for the action of $A^*$ on the boundary $\partial \mathcal{H}$, that is $A^*(\omega)\cdot \alpha(\omega)=\alpha(T\omega)$ for every $\omega \in \Omega$.\\   

We will study the {\it induced cocycle on the space of horospheres centered at $\alpha(\omega)$ for every $\omega\in \Omega$}. More precisely, for a fixed base point $p\in \mathcal{H}$ we consider the unique geodesic $\eta_{\omega}:\mathbb{R}\to \mathcal{H}$ such that $\eta_{\omega}(0)=p$ and $\eta_{\omega}(\infty)=\alpha(\omega)$ for every $\omega\in \Omega$. We define the following skew product map

\begin{eqnarray*}
V: \Omega\times \mathbb{R}&\to& \Omega\times \mathbb{R},\\
(\omega, t)&\mapsto& \left(T\omega, -b_{p, \alpha(T\omega)}(A(\omega)\cdot \eta_{\omega}(t))\right).
\end{eqnarray*}
Notice that  Lemma 4.13 implies $t=-b_{p, \alpha(\omega)}(\eta_{\omega}(t))$, hence we can see the above skew product with real fiber as a skew product acting on the space consisting into the union indexed by $\omega\in \Omega$, of the horospheres centered at $\alpha(\omega)$.\\
Let's study in more details this skew product. Since $A(\omega)$ is an isometry, the curve $A(\omega)\cdot \eta_{\omega}$ is a geodesic, and the invariance of $\alpha$ under $A^*$ implies 
\[
\big(A(\omega)\cdot \eta_{\omega}\big)(\infty)=\alpha(T\omega).
\]
Lemma 4.13 imples 
\begin{eqnarray*}
-b_{p, \alpha(T\omega)}\left(A(\omega)\cdot \eta_{\omega}(t)\right)&=&t-b_{p, \alpha(T\omega)}\left(A(\omega)\cdot \eta_{\omega}(0)\right)\\
&=&t-b_{p, \alpha(T\omega)}\left(A(\omega)\cdot p\right).
\end{eqnarray*}
 In the precedent notation, the skew product $V$ is a continuous cocycle by translations of the real line taking the form
\[
V(\omega, t)=(T\omega, t+\phi(\omega)), 
\]
where
\begin{equation}\label{phi}
    \phi(\omega)=-b_{p, \alpha(T\omega)}\left(A(\omega)\cdot p\right).
\end{equation}

\noindent{\bf Proposition 5.11. } Let $\alpha:\Omega\to \partial \mathcal{H}$ be an $A^*$ invariant section. Under the previous notation, if $(T,A)$ is $\tau$-H\"older and $\alpha$  is $\tau$-H\"older$-$Busemann  then $\phi$ is $\tau$-H\"older. Moreover if $(T,A)$ verifies the $(PPO)$ condition  then $\phi$ verifies the $(PPO)$ condition.

\begin{proof} 
Let $\omega_1, \omega_2\in \Omega$. Recall that $\alpha(T\omega_1)=A^*(\omega_1)\cdot \alpha(\omega_1)$ and $\alpha(T\omega_2)=A^*(\omega_2)\cdot \alpha(\omega_2)$. Hence, Lemma 4.14. implies that
\begin{eqnarray*}
\phi(\omega_2)-\phi (\omega_2)&=&b_{p, \alpha(T\omega_1)}(A(\omega_1)\cdot p)-b_{p, \alpha(T\omega_2)}(A(\omega_2)\cdot p)\\
&=&b_{p, \alpha(\omega_2)}(A^{-1}(\omega_2)\cdot p)-b_{p, \alpha(\omega_1)}(A^{-1}(\omega_1)\cdot p)\\
&=&b_{p, \alpha(\omega_2)}(A^{-1}(\omega_2)\cdot p)-b_{p, \alpha(\omega_2)}(A^{-1}(\omega_1)\cdot p)\\&&\quad +b_{p, \alpha(\omega_2)}(A^{-1}(\omega_1)\cdot p)-b_{p, \alpha(\omega_1)}(A^{-1}(\omega_1)\cdot p).
\end{eqnarray*}
Using Lemma 4.15., Lemmas 5.8. and Lemma 5.10.,  we obtain
\begin{eqnarray*}
|\phi(\omega_2)-\phi (\omega_1)|&\leq&|b_{p, \alpha(\omega_2)}(A^{-1}(\omega_2)\cdot p)-b_{p, \alpha(\omega_2)}(A^{-1}(\omega_1)\cdot p)|\\&&\quad +|b_{p, \alpha(\omega_2)}(A^{-1}(\omega_1)\cdot p)-b_{p, \alpha(\omega_1)}(A^{-1}(\omega_1)\cdot p)|\\ 
&\leq& d_{\mathcal{H}}(A^{-1}(\omega_2)\cdot p, A^{-1}(\omega_1)\cdot p)+Cd_{\Omega}(\omega_1, \omega_2)^{\tau}\\
&\leq& 2Cd_{\Omega}(\omega_1, \omega_2)^{\tau}.
\end{eqnarray*}
Let's verify the claim about the PPO condition. For that we will show that $V$ is a factor of $F$. More precisely we define $P(\omega, h)=(\omega, -b_{p, \alpha(\omega)}(h))$. We claim that 
\begin{equation}\label{factor}
P\circ F=V\circ P.
\end{equation}
We will give a proof of (\ref{factor}) below. Notice that (\ref{factor}) implies $P\circ F^n=V^n\circ P$. The $(PPO')$ condition $F^{n^*}(\omega^{*}, h)=(\omega^{*}, h)$ gives
\begin{eqnarray*}
P(\omega^{*}, h)&=&V^{n^*}\left(P(\omega^{*}, h)\right)\\
&=&V^{n^*}(\omega^{*}, -b_{p, \alpha(\omega^*)}(h))\\
&=&\left(\omega^{*}, -b_{p, \alpha(\omega^{*})}(h)+\sum_{j=0}^{n^*-1}\phi(T^j\omega^{*})\right).
\end{eqnarray*}
The last inequality implies $\sum_{j=0}^{n^*-1}\phi(T^j\omega^{*})=0$.\\

Let's finish the proof of this proposition by showing (\ref{factor}). Indeed, since $A$ is an isometry and $A^*(\omega)\cdot \alpha(\omega)=\alpha(T\omega)$, for every $h \in \mathcal{H} $ we have
\begin{equation}\label{factor1}
b_{A(\omega)\cdot p, \alpha(T\omega)}(A(\omega)\cdot h)=b_{p, \alpha(\omega)}(h).
\end{equation}
Moreover, Lemma 4.14. implies
\begin{equation}\label{factor2}
b_{A(\omega)\cdot p, \alpha(T\omega)}(A(\omega)\cdot h)=b_{p, \alpha(T\omega)}(A(\omega)\cdot h)-b_{p, \alpha(T\omega)}(A(\omega)\cdot p).    
    \end{equation}
The equality of  (\ref{factor1}) and (\ref{factor2}) gives
\[
-b_{p,\alpha(T\omega)}(A(\omega)\cdot h)=-b_{p, \alpha(\omega)}(h)-b_{p, \alpha(T\omega)}(A(\omega)\cdot p),
\]
which is the $\mathcal{H}$-coordinate of the equality $P\circ F=V\circ P$.
\end{proof}

\begin{proof}{\em of  Theorem 5.3. } Let $(\omega_0, h_0)\in \Omega\times \mathcal{H}$. Take $\alpha_0\in  \partial \mathcal{H}$. Let $\alpha:\Omega\to \partial\mathcal{H}$ be the $\tau$-H\"older-Busemann section that is invariant under $A^*$ and such that $\alpha(\omega_0)=\alpha_0$. The discussion in the previous section, and in particular Proposition 5.11., allows us to apply the Classic Li\v vsic Theorem  to the real valued cocycle $(T, V)$, induced by $(T, A)$  and the invariant section $\alpha:\Omega\to \partial \mathcal{H}$ on the space of horospheres.  This gives us a $\tau$-H\"older solution $u:\Omega \to \mathbb{R}$ to the cohomological equation
\begin{equation}\label{Livsic_u}
u(T\omega)-u(\omega)=\phi(\omega), 
\end{equation}
where $\phi(\omega)$ is defined at (\ref{phi}). We can also choose $u$ such that $u(\omega_0)=-b_{p, \alpha(\omega_0)}(h_0)$. By identifying $u(\omega)$ with the horosphere 
\[\mathcal{B}_{\omega}=b_{p, \alpha(\omega)}^{-1}(-u(\omega)) \ \subset \ \mathcal{H},
\]
and the fact that $V\circ P=P\circ F$, we see that the section $\omega \mapsto \mathcal{B}_{\omega}$ is invariant under the action of $F$, that is $A(\omega)\cdot \mathcal{B}_{\omega}=\mathcal{B}_{T\omega}$. Also notice that $h_0\in \mathcal{B}_{\omega_0}$.

 Let $\gamma_{\omega_0}:\mathbb{R}\to \mathcal{H}$ be a geodesic such that $\gamma_{\omega_0}(u(\omega_0))=h_0$  and $\gamma_{\omega_0}(\infty)=\alpha(\omega_0)\in \partial \mathcal{H}$. Notice that $b_{p, \alpha(\omega_0)}(\gamma_{\omega_0}(0))=0$.
 
 Since $(T, A^*)$ satisfy $(P3)$, there exists a continuous section $\beta:\Omega\to \partial{\mathcal{H}}$, that is $A^*$-invariant and such that $\beta(\omega_0)=\gamma_{\omega_0}(-\infty)$. For every $\omega\in \Omega$ let's consider  $\gamma_{\omega}:\mathbb{R}\to \mathcal{H}$ to be the unique geodesic such that
 \[
 \gamma_{\omega}(-\infty)=\beta(\omega), \quad \gamma_{\omega}(0)\in b_{p, \alpha(\omega)}^{-1}(0), \quad \gamma_{\omega}(+\infty)=\alpha(\omega). 
 \]
Due to $\gamma_{\omega}(+\infty)=\alpha(\omega)$ we obtain that for every $s\in \mathbb{R}$ we have $b_{p, \alpha(\omega)}(\gamma_{\omega}(s)))=-s$, and hence $\gamma_{\omega}(u(\omega))\in \mathcal{B}_{\omega}$. Moreover, this also yields $\{\gamma_{\omega}(u(\omega))\}=\gamma_{\omega}\cap \mathcal{B}_{\omega}$.

As $A$ is an isometry, and $A^*(\omega)\cdot \beta(\omega))=\beta(T\omega)$, $A^*(\omega)\cdot \alpha(\omega)=\alpha(T\omega)$, we conclude that 
\[
A(\omega)\cdot \gamma_{\omega}=\gamma_{T\omega}, \quad A(\omega)\cdot \gamma_{\omega}(u(\omega))\in \gamma_{T\omega}.
\]
The previous line, togheter with $A(\omega)\cdot \mathcal{B}_{\omega}=\mathcal{B}_{T\omega}$ allows to conclude that \[
A(\omega)\cdot \gamma_{\omega}(u(\omega))\in \mathcal{B}_{T\omega}\cap \gamma_{T\omega}.
\]
Hence, $A(\omega)\cdot \gamma_{\omega}(u(\omega))=\gamma_{T\omega}(u(T\omega))$. It is time to  define 
\[
s_{\omega_0, h_0}(\omega)=\gamma_{\omega}(u(\omega)), 
\]
which is a  $A$-invariant section, that verifies  $s_{\omega_0, h_0}(\omega_0)=\gamma_{\omega_0}(u(\omega_0))=h_0$, 
that is, it passes through $(\omega_0, h_0)$. \\

\noindent{\it Claim.} As defined before, $\omega \mapsto s_{\omega_0, h_0}(\omega)$ is continuous. \\

Since $\omega \mapsto u(\omega)$ is continuous then it suffices to show that $\omega\mapsto \gamma_{\omega}(0)$ is continuous. We will use twice the Lemma 4.15. Recall that $\gamma_{\omega}(0)$ is the unique point on $\gamma_{\omega}$ such that $b_{p, \alpha(\omega)}(\gamma_{\omega}(0))=0$. Given $\overline{\omega}\in \Omega, \varepsilon>0$,  let's apply Lemma 4.10. to the geodesic $\gamma_{\overline{\omega}}$, the point $\gamma_{\overline{\omega}}(\varepsilon)$ (respectively $\gamma_{\overline{\omega}}(-\varepsilon)$) and the radius $r=\varepsilon/2$. By hypothesis,  for fixed $h$ the Busemann functions $\omega\mapsto b_{p, \alpha(\omega)}(h)$ are uniformly H\"older, for $h$ in a bounded set. Hence for $\omega$ close enough to $\overline{\omega}$ we have 
\begin{eqnarray*}
b_{p, \alpha(\omega)}\big |_{B(\gamma_{\overline{\omega}}(-\varepsilon), \varepsilon/2)}&>&b_{p, \alpha(\overline{\omega})}(-\varepsilon)/2=\varepsilon/2>0, \\
b_{p, \alpha(\omega)}\big |_{B(\gamma_{\overline{\omega}}(\varepsilon), \varepsilon/2)}&<&b_{p, \alpha(\overline{\omega})}(\varepsilon)/2=-\varepsilon/2<0.
\end{eqnarray*}
Henceforth, for $\omega$ close enough to $\overline{\omega}$, the geodesic $\gamma_{\omega}$ traverses through the small balls $B(\gamma_{\overline{\omega}}(\varepsilon),\varepsilon)$ and  $ B(\gamma_{\overline{\omega}}(-\varepsilon), \varepsilon)$. Hence, the  function 
\[
t\mapsto b_{p, \alpha(\omega)}(\gamma_{\omega}(t))
\]
reaches its unique zero inside the ball $B(\gamma_{\overline{\omega}}(0), 3\varepsilon)$, thanks to the convexity of balls. This completes the claim.\\

In order to finish the proof of Theorem 5.3., we need to show that for every $\omega_0, \omega\in \Omega$ the map $ h\mapsto s_{\omega_0, h}(\omega)$ is continuous. Denote $u_{h}$ the solution of (\ref{Livsic_u}) for $h_0=h$. Recall that $u_h(\omega_0)=-b_{p, \alpha(\omega_0)}(h)$. Also, from section \ref{CLT} we know that for distinct $h, \overline{h}\in \mathcal{H}$ the two solutions to (\ref{Livsic_u}) differ by a constant, hence
\begin{eqnarray}
        u_{h}&=&u_{\overline{h}}+\left(u_h(\omega_0)-u_{\overline{h}}(\omega_0)\right),\\
        &=&u_{\overline{h}}+\left(b_{p, \alpha(\omega_0)}(\overline{h})-b_{p, \alpha(\omega_0)}(h)\right).
\end{eqnarray}
Thus we obtain that $h\mapsto u_h$ is uniformly continuous. The remaining part of the proof relies just in geometric arguments, whose continuity closely follows the uniqueness of geodesics with endpoints in $\mathcal{H}\cup \partial\mathcal{H}$. Let's recall how the construction of $s_{\omega_0, h}(\omega)$ depends on $h$. We construct a geodesic $\gamma_{\omega_0}:=\gamma_{\omega_0, h}$ that joins $h\in \mathcal{H}$ to $\alpha(\omega_0)\in \partial \mathcal{H}$. We denote by $\beta_0:=\beta_{0,h}$ the point $\gamma_{\omega_0, h}(-\infty)$. Lemma 4.15. says that $h\mapsto \beta_{0,h}$ is continuous. Afterwards, we consider the continuous $A^*$-invariant section $\beta_{h}:\Omega\to \partial \mathcal{H}$ that passes through $\beta_{0, h}$. This section depends (pointwise on $\omega$) continuously on $h$. For each $\omega\in \Omega$ we consider the geodesic $\gamma_\omega:=\gamma_{\omega, h}$ that joins $\beta_{h}(\omega )$ with $\alpha(\omega)$. Continuity (with respect to $h$) of $\gamma_{\omega, h}(u_h(\omega))$ follows.
\end{proof}

\section*{ Appendix. A metric on the Gromov boundary, strong hyperbolicity and considerations on the H\"older conditions}
In this appendix we review a class of  metric spaces where it is possible to establish a metric on the boundary at infinity, so that the H\"older conditions on the solutions on the boundary of the cohomnological equation can be expressed in a simpler way. A general reference for strong hyperbolic spaces is \cite{NIKA}.\\

\noindent{\bf Definition A1. }
Let $\mathcal{H}$ be a geodesic  space. Given three points $x,y,p \in \mathcal{H}$ the \textit{Gromov Product} is defined as
\begin{equation}\label{Gromov_product}
(x,y)_p=\frac{1}{2}(d_{\mathcal{H}}(x,p)+d_{\mathcal{H}}(y,p)-d_{\mathcal{H}}(x,y)).
\end{equation}

In some sense, the number $(x, y)_p$ measures how thin is the traingle $xyp$. We can define another boundary for a Gromov hyperbolic space using the the Gromov product. Let be $\mathcal{H}$ be a Gromov hyperbolic space and a base point $p\in \mathcal{H}$. We say that a sequence $(x_i)\subset \mathcal{H}$ converges to infinity if $(x_i,x_j)_p\to \infty$ when $i,j\rightarrow \infty$. Two such sequences $x_i,y_j\in \mathcal{H}$ are related if $(x_i,y_j)_p\rightarrow \infty$. This defines an equivalence relation among sequences that converge to infinity (transitivity is a consequence of $\delta$-hyperbolity). This defines a  boundary $\partial_s \mathcal{H}$ as the set of equivalences classes for this relation.  The map $(x_n)\rightarrow \lim_n x_n $ is a bijection between $\partial_s \mathcal{H}$ and $\partial \mathcal{H}$. Indeed,  a sequence $(x_i) \subset \mathcal{H}$  converges to a point in $\partial \mathcal{H}$ if and only if converges to infinity.\\

Now we extend $(\cdot,\cdot)_p$ to the boundary. We could hope that the extension can be defined simply as $(x,y)_p= \lim_{i,j\rightarrow \infty}(x_i,y_j)_p$. Nevertheless,  it isn't possible in general. The limit could not exist or even depend on the representative sequences  (see \cite{BRIDSON}, 3.16). Notwithstanding, there is a large list of spaces where this limit does holds, and in that list we can find familiar examples, like the traditional hyperbolic space and trees.\\

In general Gromov hyperbolic spaces, one can define the Gromov product for points  $x, y \in \partial \mathcal{H}$, by
\[
(x,y)_p=\sup \liminf_{i,j\rightarrow \infty}(x_i,y_j)_p,
\]
where $p\in \mathcal{H}$ is a fixed base point and the supremum is taken over all sequences $x_i,y_j \in \mathcal{H}$ such that $\lim_i x_i=x$ and $\lim_j y_j=y$. 
Recall that given a triangle with vertices $x,y,z\in \mathcal{H}$, the more large is $(x,y)_z$ more thin is the triangle. With this in mind, we define now a topology for the boundary.  \\

\noindent{\bf Definition A2.} A sequence $\xi_n\in \partial H$ converges to a point $\xi\in \partial H$ if 
\[
(\xi_n,\xi)_p\rightarrow \infty.
\]
This definition doesn't depend on the base point $p$. The next step is to construct a metric for the boundary. Following  \cite{BRIDSON} section 3.19, for $\varepsilon>0$ we define
\begin{equation}\label{metric_boundary}
\varrho_{\varepsilon}(x,y)=e^{-\varepsilon(x,y)_p}.
\end{equation}

\noindent{\bf Definition A3. }
A hyperbolic metric space   $\mathcal{H}$ is said to be {\it strongly hyperbolic} of parameter $\varepsilon>0$ if for all points $x,y,z\in \mathcal{H}$ and a base point $p\in \mathcal{H}$ the following inequality holds
\[
\varrho_{\varepsilon}(x,y)\leq \varrho_{\varepsilon}(x,z)+\varrho_{\varepsilon}(y,z).
\]
The recent article  \cite{das2017geometry}  provides a lot of examples of strong hyperbolic spaces. These are general enough as the following  shows.\\

\noindent{\bf Theorem A4., see \cite{das2017geometry}.}
Every CAT$(-1)$ space is strongly hyperbolic.\\

Since every $\mathbb{R}$-tree is $\textup{CAT}(-1)$, they are strongly hyperbolic. The Poincar\'e plane $\mathbb{H}^n$ is strongly hyperbolic.\\

\noindent{\bf Theorem A5., see \cite{NIKA} Th. 4.2. } Let $\mathcal{H}$ be a strongly hyperbolic space and $p$ a basepoint in $\mathcal{H}$. Then the Gromov product $(\cdot,\cdot)_p$ extends continuously to the boundary $\partial \mathcal{H}$. Moreover,  for every $x,y\in\partial \mathcal{H}$ and sequences $x_i\rightarrow x$, $y_i\rightarrow y$ 
 \begin{equation}\label{limit}
 (x,y)_p=\lim_{i\rightarrow \infty}(x_i,y_i)_p.
 \end{equation}
 Further, the definition doesn't depend on the representative sequences. Moreover, $\mathcal{H}$ is Gromov hyperbolic and $\varrho_{\varepsilon}$ (as defined in (\ref{metric_boundary}), is a metric on $\partial\mathcal{H}$.\\

Another consequence is that Busseman  functions can be rewritten only in therms of the Gromov product. \\

\noindent{\bf Lemma A6.} In the current notation the following holds.
\begin{equation}\label{buseman_strong}
b_{h,\xi}(p)=2(\xi,h)_{p}-d(h,p).
\end{equation}

\noindent{\bf Lemma A7.}
If $\mathcal{H}$ is a strongly hyperbolic space and $g\colon \mathcal{H}\rightarrow \mathcal{H}$ is an isometry, the extended action $g^*:\partial \mathcal{H}\rightarrow \partial \mathcal{H}$ is a Lipsichtz map.\\

\noindent{\bf H\"older conditions.} Let's explore the key estimate given by Proposition 5.11. in order to deduce that the cocycle by translations (\ref{phi}) is H\"older. Lemma 4.14. implies
\begin{eqnarray*}
\phi(\omega_2)-\phi(\omega_1)&=& b_{A^{-1}(\omega_2)\cdot p, \alpha(\omega_2)}(p)-b_{A^{-1}(\omega_1)\cdot p, \alpha(\omega_1)}(p),\\
&=&2(\alpha(\omega_2), A^{-1}(\omega_2)\cdot p,)_p-d_{\mathcal{H}}(A{^-1}(\omega_2)\cdot p, p)+\\
&&- 2(\alpha(\omega_1), A^{-1}(\omega_1)\cdot p,)_p+d_{\mathcal{H}}(A^{-1}(\omega_1)\cdot p, p).
\end{eqnarray*}
Using (\ref{buseman_strong}) we can estimate
\begin{eqnarray*}
|\phi(\omega_2)-\phi(\omega_1)|&\leq&|d_{\mathcal{H}}(p, A^{-1}(\omega_2)\cdot p)-d_{\mathcal{H}}(p, A^{-1}(\omega_1)\cdot p)|+\\ 
&&+2\big\{|(\alpha(\omega_2), A^{-1}(\omega_2)\cdot p)_p-(\alpha(\omega_2),A^{-1}(\omega_1)\cdot p)_p|+\\ 
&&+|(\alpha(\omega_1), A^{-1}(\omega_1)\cdot p)_p-(\alpha(\omega_2),A^{-1}(\omega_1)\cdot p)_p |\big \}.
\end{eqnarray*}
Following (\ref{Gromov_product}) and (\ref{limit}), we know that for $\xi\in \partial{\mathcal{H}}$ and $x, y\in \mathcal{H}$ one has $|(\xi, x)_p-(\xi, y)_p|\leq d_{\mathcal{H}}(x, y)$. Hence we have
\begin{eqnarray*}
|\phi(\omega_2)-\phi(\omega_1)|&\leq&3d_{\mathcal{H}}(A^{-1}(\omega_2)\cdot p, A^{-1}(\omega_1)\cdot p)+\\ 
&&+2|(\alpha(\omega_1), A^{-1}(\omega_1)\cdot p)_p-(\alpha(\omega_2),A^{-1}(\omega_1)\cdot p)_p |.
\end{eqnarray*}
We can use (\ref{buseman_strong}) in order to obtain that \begin{eqnarray*}
\omega_1,\omega_2&\longmapsto& (\alpha(\omega_2), A^{-1}(\omega_1)\cdot p)_p,\\
\omega_1&\longmapsto& (\alpha(\omega_1), A^{-1}(\omega_1)\cdot p)_p
\end{eqnarray*}
are uniformly bounded (from above). Hence, 
\begin{eqnarray*}
\omega_1,\omega_2&\longmapsto& \varrho_{\varepsilon}(\alpha(\omega_2), A^{-1}(\omega_1)\cdot p),\\
\omega_1&\longmapsto& \varrho_{\varepsilon}(\alpha(\omega_1), A^{-1}(\omega_1)\cdot p)
\end{eqnarray*}
are  uniformly  away from zero. Then there exists $C>0$ such that
\begin{eqnarray*}
|(\alpha(\omega_1), A^{-1}(\omega_1)\cdot p)_p-(\alpha(\omega_2),A^{-1}(\omega_1)\cdot p)_p |&=&\\
\frac{1}{\varepsilon}|\log \varrho_{\varepsilon}(\alpha(\omega_1), A^{-1}(\omega_1)\cdot p)-\log \varrho_{\varepsilon}(\alpha(\omega_2), A^{-1}(\omega_1)\cdot p)|&\leq&\\
\frac{C}{\varepsilon}| \varrho_{\varepsilon}(\alpha(\omega_1), A^{-1}(\omega_1)\cdot p)- \varrho_{\varepsilon}(\alpha(\omega_2), A^{-1}(\omega_1)\cdot p)|&\leq&\\
\frac{C}{\varepsilon}\varrho_{\varepsilon}(\alpha(\omega_1), \alpha(\omega_2)).
\end{eqnarray*}
Summarizing we obtain that
\begin{eqnarray*}
|\phi(\omega_2)-\phi(\omega_1)|&\leq&3d_{\mathcal{H}}(A^{-1}(\omega_2)\cdot p, A^{-1}(\omega_1)\cdot p)+\\ 
&&+\frac{2C}{\varepsilon}\varrho_{\varepsilon}(\alpha(\omega_1), \alpha(\omega_2)).
\end{eqnarray*}
Since the core of the proof of the Theorem 5.3. is the fact that $\phi$ is H\"older (in order to use the classic Li\v vsic Theorem), we can state the next Proposition in the context of strongly hyperbolic metric spaces. It implies in particular that Theorem 5.3. holds if we replace the H\"older-Busseman condition on the boundary by a classic H\"older condition in the  $\varrho_{\varepsilon}$ metric on the boundary.\\

\noindent{\bf Proposition A8.} Let $\mathcal{H}$ be a strongly hyperbolic metric space. Let $\alpha:\Omega\to \partial \mathcal{H}$ be an $A^*$ invariant section. Under the previous notation, if $(T,A)$ is $\tau$-H\"older and $\alpha$  is $\tau$-H\"older (for the metric $\varrho_{\varepsilon}$), then $\phi$ is $\tau$-H\"older.\\

%\noindent{\it Acknowledgments.} We appreciate the helpful comments of Jairo Bochi and Alejandro Kocsard. 


\begin{thebibliography}{}


\bibitem{AVKOLI}{\sc A. Avila,  A. Kocsard and X. Liu.} Li\v vsic theorem for diffeomorphism cocycles. {\em Geometric and Functional Analysis (GAFA)}  {\bf 28} (2018), 943-964. 

\bibitem{AVKR06}{\sc A. Avila \& R. Krikorian.} Reducibility or nonuniform 
hyperbolicity for quasiperiodic Schr\"{o}dinger cocycles. {\em Ann. of Math.} 
{\bf 164} (2006), 911-940.

%\bibitem{ballmann} {\sc A. Ballmann.} {\em Lectures on Spaces of 
%Nonpositive Curvature.} DMV Seminar {\bf 25}, Birkh\"auser (1995).

%\bibitem{BOST} {\sc J.-B. Bost.} Tores invariants des syst\`emes dynamiques hamiltoniens (d'apr\`es Kolmogorov, Arnol'd, Moser, R\"ussmann, Zehnder, Herman, P\"oschel,...), S\'eminaire Bourbaki, vol. 1984/85, {\em Ast\'erisque} 133-134 (1986), p. 113-157.

\bibitem{BRIDSON} {\sc Bridson, M.R. and H{\"a}fliger, A.}
{\em Metric Spaces of Non-Positive Curvature}, Grundlehren der mathematischen Wissenschaften, Springer Berlin Heidelberg (2011).


\bibitem{CONAPO} {\sc D. Coronel, A. Navas and M. Ponce.} On bounded  cocycles  isometries over a minimal dynamics. {\em  Journal of Modern Dynamics}, vol {\bf 7}, pages: 45 - 74, Issue 1, March 2013.

\bibitem{das2017geometry} {\sc Das, T. and Simmons, D. and Urba{\'n}ski, M.},
{\em Geometry and Dynamics in Gromov Hyperbolic Metric Spaces. } {Mathematical Surveys and Monographs},
 {American Mathematical Society, 2017.}


%\bibitem{dela86} {\sc R. de la Llave, J. M. Marco and R Moriy\'on.} Canonical perturbation theory of Anosov Systems and regularity for the Liv\v sic cohomology equation, {\em Ann. of Math. } (2) {\bf 123}(3) (1986), 537-611.

\bibitem{delallave2010} {\sc R. de la Llave \& A. Windsor.} Liv\v sic theorems for non-commutative
groups including  diffeomorphism groups and results on the existence of  conformal structures for
Anosov systems. {\em Ergodic Theory and Dynam. Systems} {\bf 30}, no. {\bf 4} (2010), 1055-110.

\bibitem{dolgo05} {\sc D. Dolgopyat.} Liv\v sic theory for compact group extensions of hyperbolic systems. {\em Mosc. Math. J.} {\bf 5}(1)(2005), 55-67.

%\bibitem {eliasson} {\sc L. H. Eliasson.} Perturbations of stable invariant tori for Hamiltonian systems. {\em Ann. Scuola Norm. Sup. Pisa Cl. Sci. }(4), 15(1):115-147 (1989), 1988.

\bibitem{EO}{\sc  P. Eberlein and B. O’Neill.} Visibility manifolds. {\em Pacific J. Math.} {\bf 46} (1973), 45–109.

\bibitem{FK} {\sc B. Fayad and R. Krikorian}. Exponential growth of product of matrices in $SL(2, \mathbb{R})$. {\em Nonlinearity}, {\bf 21}(2) (2008), 319-323.

\bibitem{forfla}{\sc G. Forni and L. Flaminio.} On the cohomological equation for nilflows. {\em J. Mod. Dyn.} {\bf 1} (1)(2007), 37-60.


\bibitem{GOHE} {\sc W. H. Gottschalk and G. A. Hedlund.} { Topological Dynamics.} 
{\em Amer. Math. Soc.}, Providence, R. I. (1955).

\bibitem{Gromov1987} {\sc Gromov, M.},
 {\em Hyperbolic Groups},
Essays in Group Theory, pp. 75-263, 
Springer New York, 1987.


\bibitem{KAL10} {\sc B. Kalinin.} Liv\v sic theorem for matrix cocycles.
{\em Ann. of Math.} {\bf 173}, no {\bf 2} (2011), 1025-1042.

\bibitem {KS}{\sc B. Kalinin \& V. Sadovskaya}. Linear cocycles over hyperbolic systems 
and criteria of conformality. {\em J. Mod. Dyn.} {\bf 4}, no. {\bf 3} (2010), 419--441.

\bibitem{KAR_MAR} {\sc Karlsson, A., Margulis, G.A.} A multiplicative ergodic theorem and nonpositively curved spaces. {\em Comm. Math. Phys. }{\bf 208}(1) (1999), 107-123.

\bibitem{KAKO} {\sc A. Katok and A. Kononenko.} Cocycles' stability for partially hyperbolic systems. {\em Math. Res. Lett.} {\bf 3} (1996), no. 2, 191-210.

 \bibitem{KARO} {\sc A. Katok.} Cocycles, cohomology and combinatorial constructions in ergodic theory. In collaboration with E. A. Robinson. {\em Jr. Proc. Sympos. Pure Math.} {\bf 69}, Smooth ergodic theory and its applications (Seattle, WA, 1999), 107-173, Amer. Math. Soc., Providence, RI, 2001.

\bibitem{KOPO}{\sc A. Kocsard and R. Potri\'e.} Li\v vsic theorem for low-dimensional diffeomorphism cocycles. {\em Commentarii Mathematici Helvetici} {\bf 91} (2016), 39-64.


\bibitem{LIV72} {\sc A.N. Liv\v sic.} Cohomology of dynamical systems. {\em Math. USSR Izvestija}
{\bf 6} (1972), 1278-1301.

\bibitem{MMY} {\sc S. Marmi, P. Moussa \& J.-C. Yoccoz.} The cohomological 
equation for {R}oth-type interval exchange maps. {\em J. Amer. Math. Soc.} 
{\bf18}, No. {\bf 4} (2005), 823-872. 


\bibitem{NavasPonce} {\sc A. Navas and M. Ponce.} A Liv\v sic type theorem for germs of analytic diffeomorphisms. {\em Nonlinearity} {\bf 26} (2013) 297-305.

\bibitem{NIKA} {B. Nica and J. Spakula.} Strong Hyperbolicity. {\em Groups Geom. Dyn.} {\bf 10} (2016) 951-964.

\bibitem{nittor1} {\sc V. Nitica and A. T\"or\"ok.} Local rigidity of certain partially hyperbolic actions of product type. {\em Ergodic Theory and Dys. Sys.} {\bf 21}(4) (2001), 1213-1237. 

%\bibitem{nittor2} {\sc V. Nitica and A. T\"or\"ok.} Regularity of the transfer map for cohomologous cocycles. {\em Ergodic Theory and Dys. Sys.} {\bf 18}(5) (1998), 1187-1209. 

\bibitem{poll99} {\sc M. Pollicott and M. Yuri.} Regularity of solutions to the measurable Liv\v sic equation. {\em Trans. Amer. Math. Soc.} {\bf 351}(2)(1999), 559-568. 



%\bibitem{poschel} {\sc J. P\"oschel.} On elliptic lower-dimensional tori in Hamiltonian systems. {\em Math. Z.} {\bf 202} (1989), 559-608.


  \bibitem{WILK}{\sc A. Wilkinson.} The cohomological equation for partially hyperbolic diffeomorphisms. {\em Asterisque} {\bf 358} (2013), 75-165.
  
%\bibitem{yoccoz2} {\sc J.-C. Yoccoz.} Some questions and remarks about $SL(2,\mathbb{R})$ 
cocycles. {\em Modern dynamical systems and applications}. Cambridge Univ. Press, 
Cambridge (2004), 447-458.

\bibitem{YOCC} {Yoccoz, J.-Ch.} {\it Interval exchange maps and translation surfaces. Homogeneous flows, moduli spaces and arithmetic}. Clay Math. Proc., 10, Amer. Math. Soc., pp. 1-69, Providence, RI, 2010.



\end{thebibliography}
\end{document}